\documentclass{article}

\begin{document}


\input{amssym.def}
\input{amssym.tex}

 \newcommand{\R}{{\Bbb R}}
 \newcommand{\Ha}{{\cal H}}
 \newcommand{\Z}{{\Bbb Z}}
 \newcommand{\C}{{\Bbb C}}
 \newcommand{\N}{{\Bbb N}}
\newcommand{\F}{{\cal F}}
\newcommand{\Q}{{\Bbb Q}}
\newcommand{\Ro}{{\cal R}}

\title{CMC--Surfaces, $\varphi$--Geodesics\\
                     and\\
     Caratheodory'sche Vermutung}

\bigskip

 \author{Igor ~Nikolaev\\
 The Fields Institute\\
222 College street, Toronto\\
 M5T 3J1 ~Canada\\
{\sf E-mail: inikolae@fields.utoronto.ca}}

 \maketitle

\newtheorem{thm}{Theorem}
\newtheorem{lem}{Lemma}
\newtheorem{dfn}{Definition}
\newtheorem{rmk}{Remark}
\newtheorem{cor}{Corollary}
\newtheorem{prp}{Proposition}

\begin{abstract}
A short proof of the Caratheodory conjecture about 
index of an isolated umbilic on the convex 2--dimensional
sphere is suggested.


\vspace{7mm}

{\it Key words and phrases:  Quadratic Differential, $\varphi$--Metric,
 Umbilical Point}

\vspace{5mm}
{\it AMS (MOS) Subj. Class.:  30F30, 53A07, 58F10.}
\end{abstract}

\section*{Introduction}
The constant mean curvature (CMC--) surfaces in $E^3$ are known to admit a continuous family
of local, non--trivial, isometric deformations preserving mean curvature
of the surface ($H$--deformations).
In the case when surface is compact Umehara \cite{Ume} showed that the converse is also true.

The maximal and minimal curvature lines of CMC--surface form an orthogonal net which is called
 {\it r\'eseau de Bonnet}, cf Cartan \cite{Car}. The Bonnet Theorem says that if the CMC--surface is
simply connected and umbilic--free, then under $H$--deformations the orthogonal net
``rotates'' through a constant angle which can be taken as a parameter of deformation.

If the CMC--surface is not simply connected or umbilic--free, Cartan seems to be the
first to ask about possible scenario of evolution of {\it r\'eseau de Bonnet} 
under the $H$--deformations.

In the present note we study
\footnote{For the reasons which will be clear later.}
evolution of the orthogonal nets in the case when CMC--surface
is simply connected with a single umbilic or, equivalently, doubly connected
and umbilic--free. Namely, if we ``pinch'' the umbilic, the CMC--surface
becomes an annulus whose points undergo $H$--deformations accordingly with
the Bonnet Theorem. In general, the rotation angle is no longer constant 
at all the points because  annulus cannot be covered by a single chart.

However, the Bonnet Theorem implies that every curve of the orthogonal net
is a $\varphi$--geodesic line whatever $H$--deformations are applied to
a CMC--surface. Metric $\varphi$ is given by the linear element 
$ds=|\varphi||dz|$, where $\varphi dz^2$ is a holomorphic quadratic
"differential" associated to the CMC--surface. Of course, $\varphi(0)=0$
at the umbilical point.

This observation is crucial,  because the $\varphi$--geodesics
near $n$--th order zero of a holomorphic quadratic form are
well--understood due to the works of Strebel \cite{Str}.
Roughly speaking, the $\varphi$--geodesics fill--up the annulus
either by  ``hyperbolas'' or  ``radii''. Therefore, possible configurations
of  {\it r\'eseau de Bonnet} near the umbilic looks like a singularity
with the finite number of hyperbolic and parabolic sectors.

Despite  independent interest, the orthogonal nets are auxiliary
for us. We postulate different fact here:  $H$--deformations of
orthogonal nets give an amazingly simple proof
to the {\it Caratheodory'sche Vermutung (Conjecture)}:
\begin{thm}\label{C}
Let $S^2$ be a $C^{\infty}$ surface which bounds a convex compact body 
in the Euclidean space $E^3$. Then $S^2$ has at least
two umbilical points. In other words, 
the Euler-Poincar\'e index of isolated umbilical point
is at most $+1$.
\end{thm}
(A short overview of this conjecture can be found in \cite{AVL};
see also \cite{Bol}, \cite{Ham}, \cite{Klo}.)

\bigskip\noindent
{\sf Acknowledgments.}  I am grateful to the referee for a careful reading of 
this note.

\section{$\varphi$--Geodesics}
Until further indications, $M$ is a simple domain of the complex
parameter $z$. Let us consider the holomorphic functions $\varphi(z)$ 
vanishing at the unique point of $M$ which we identify with $0$.
An  order $n\ge 1$ is assigned to $0$, if there exists a complex
constant $a\ne 0$ such that $\varphi(z)=az^n+O(|z|^{n+1})$.

Flat metric $\varphi$ with the cone singularity of angle $(n+2)\pi$ 
is given by the formula
\displaymath
|ds|=|\varphi||dz|,
\enddisplaymath
provided $\varphi(z)dz^2$ is a quadratic form on $M$.
By a {\it $\varphi$--geodesic line} in $M$ one understands 
the line conisting of the shortest arcs relatively metric $\varphi$.
Any two points in $M$ (including $0$) may be joined by the unique
$\varphi$--geodesic line. Strebel classified the possible types
of $\varphi$--geodesics in the neighborhood of $n$--th order zero
by proving the following lemma.
\begin{lem}(\cite{Str})\label{geo}
Any two points in a neighborhood $M$ of $n$--th order zero
of holomorphic 2--form $\varphi(z)dz^2$ can be joined by a
unique $\varphi$--geodesic. Moreover, each $\varphi$--geodesic
is either an arc defined by the equation $Arg~\varphi(z)dz^2=Const$,
or is composed of the two radii centered in $0$ with the minimal
angle $\ge 2\pi/(n+2)$.
\end{lem}
The foliation $\F$ on $M\backslash 0$ is said to be {\it geodesic}
if every leaf of $\F$ is a $\varphi$--geodesic line. Before we
state the general lemma on the structure of geodesic foliations,
let us consider an example when all $\F$'s can be obtained
by a "brute force".

\bigskip\noindent
{\scriptsize
If $0$ is a double zero, then the $\varphi$--metric is given by
the linear element $ds^2=(u^2+v^2)(du^2+dv^2)$ where $u+iv$ is a natural
parameter. The metric $|ds|$ is Liouville's
and the geodesic lines in this metric are completely integrable.
The general integral is known to be of the form
\displaymath
\int{du\over\sqrt{u^2-a}} \pm\int{dv\over\sqrt{v^2+a}}=a',
\enddisplaymath
where $a, a'$ are two independent constants. Easy calculations
show that no information will be lost if we suppose $a=0$. The
integral takes the form $\ln|u|\pm\ln|v|=a'$. The geodesic
foliation is described by two "families of curves": $v=Cu$ and 
$v=C/u$, where C is an arbitrary constant. Thus, $\F$
near a double zero is either the "node" with the geodesics radii tending to
$0$, or the "saddle" with four sectors filled--up by the geodesic "hyperbolas".
}

\bigskip\noindent
Let $w$ be a finite "word" on the alphabet consisting of two symbols $h$ and $p$.
We introduce the elementary operations on $w$:

\medskip\noindent
(i) a cyclic permutation of the symbols in $w$, and

\smallskip\noindent
(ii) a contraction of the $p$--symbol: $p^2=p$.

\bigskip\noindent
Two words are {\it equivalent} $w_1\sim w_2$ if and only if
$w_2$ can be obtained from $w_1$ by the elementary operations. 
The equivalence class of word $w$ is denoted by $[w]$.

Fix an integer number $n\ge1$. To every symbol $h$ in $w$ we 
assign a {\it weight} $|h|=2\pi/(n+2)$. To every symbol $p$
we assign the weight $|p|=\alpha_i$, where $\alpha_i$
is a positive real. The weight of $w$ is an additive
function equal to the sum of weights of the symbols entering $w$.
The equivalence class $[w]$ is called {\it normalized} if $|w|=2\pi$
for all $w\in[w]$. (Note that the weight of $w$ is one and the same
for all $w\in[w]$.) 
\begin{lem}\label{fol}
Let $h$ and $p$ stay for the hyperbolic and the parabolic sectors
of the singularity $w$, respectively. 
We encode the singular point $w$ by a sequence of symbols $h$ and $p$
in the order the $h$-- and the $p$--sectors occur when turning clockwise
around the singularity. Then:

\medskip\noindent
(i) each $\varphi$--geodesic foliation $\F$ is topologically
equivalent to the singularity $w$ of a normalized equivalence class $[w]$;

\smallskip\noindent
(ii) each normalized equivalence class $[w]$ can be realized as
a $\varphi$--geodesic foliation $\F$ with the singularity $w\in [w]$
in a neighborhood of $n$--th order zero of $\varphi$ for some
$n\ge 1$.

\end{lem}
{\it Proof.}
Denote by $M$ a neighborhood of the $n$--th order zero of $\varphi$. Let us
introduce a partial order for the points $x,y\in M$: $x\le y$ if and only if
$Arg~x\le Arg~y$. 
If $x\in M$ is an arbitrary point, then by Lemma \ref{geo} the  $\varphi$--geodesic
line through $x$ is either (i) the hyperbola $Arg~\varphi dz^2=~Const$ or 
(ii) the radius $Ox$.
Let us consider the first possibility.

(i) The hyperbola $Arg~\varphi dz^2=~Const$ must tend to the asymptotic rays $Oz_1,
Oz_2$ with $z_1<x<z_2$, enclosing the angle $2\pi/(n+2)$. Clearly, the 
only possibility to the geodesic foliation $\F$ is to form a hyperbolic
sector $z_1Oz_2$. Of course, along $Oz_1$ and $Oz_2$ $Arg~\varphi dz^2$ is constant.

(ii) Let $Ox$ be the geodesic radius through $x$, distinct from the boundary
radii of the hyperbolic sector. Then through the nearby points $|x-y|<\varepsilon$ 
one can draw the geodesic radii $Oy$'s. Denote by $y_1Oy_2$ the maximal connected parabolic 
sector filled-up with the geodesic radii. Clearly, $y_1< x< y_2$.
The angle enclosed between two boundary radii, we denote
by $\alpha$. In general, $0\le\alpha\le 2\pi$.

If the hyperbolic sector $h$ is followed by another hyperbolic sector $h$, 
we write this as $hh$. If $h$ is
 followed by a parabolic sector, we put it as $hp$. A parabolic sector
 $p$ followed by the parabolic sector $p$, gives a larger parabolic sector $p=pp$
 and the contraction rule (ii) follows. Of course, the "weights" of the sectors
 are equal to the angles swept by the sectors.

Finally, according to the definition of normalized equivalence class,
each singularity consists of sequence of parabolic and hyperbolic
sectors; every curve in these sectors is a geodesic arc.

The part (ii) of Lemma \ref{fol} is proved by the similar argument.   
$\square$

\section{CMC--Surfaces}
Every smooth immersion $f:M\to E^3$ of an orientable surface $M$
into the Euclidean space $E^3$ induces a Riemann structure
on $M$; let $z=u+iv$ be the corresponding local parameter.
With respect to $z$ the first
fundamental form can be written as $ds^2=e^{2\lambda}|dz|^2$.

If $ldu^2+2mdudv+ndv^2$ is the second fundamental form, we
consider a complex quadratic form $\varphi dz^2$, such that
$\varphi(z)={1\over 2}(l-n)-im$. The Mainardi--Codazzi equations imply
that $\varphi$ is holomorphic on $M$ if and only if $f(M)$
is a CMC--surface. Locally, along the lines of minimal and
maximal curvature $Arg~\varphi dz^2=0$ and $\varphi(0)=0$
at  the umbilic points.

A {\it continuous deformation} $f_t$ of the immersion $f=f_0$
is the isometry of surface $M$ such that $M\times[0,1]\to E^3$
is a continuous mapping. The continuous deformation $f_t$ is called
an {\it $H$--deformation} if $H_t=H$ for all $t\in[0,1]$,
where $H:M\to\R$ is the mean curvature function.

The CMC--surfaces are known to admit a non--trivial $H$--deformations
and in the case of compact surfaces, they are the only ones
with such a property. Of course, there are known many examples
of compact CMC--surfaces of genus $g>0$.

What happens with the lines of principal (i.e. minimal or
maximal) curvature of the CMC--surface during an $H$--deformation?
If $M$ is a local CMC--surface without umbilics, the principal
curvature lines of $f_0(M)$ and $f_t(M)$ form two families
of the parallel lines intersecting each other with the constant
angle proportional to the parameter $t$ (the Bonnet Theorem, see
e.g. \cite{Che}). Note, that if we fix the $\varphi$--metric
on $M$ corresponding to $f_0(M)$, then the principal
curvature lines of $f_t(M)$ coincide with the $\varphi$--geodesic lines
of the inclination $t$. If the umbilical points are allowed, then a
law is given by the following lemma.
\begin{lem}\label{angle}
Suppose that $M_0=f_0(M)$ is a canonical CMC--surface with the quadratic
function $\varphi=z^n, n\ge1$. Let $\varphi$ be a metric on $M$ corresponding
to $M_0$. If $M_t=f_t(M)$ is an
$H$--deformation of $M_0$, then one of the two principal curvature
lines of $M_t$ coincide with the $\varphi$--geodesic lines on $M$
for any $t\ge 0$.
\end{lem}
{\it Proof.} In the polar coordinate system the coefficients 
of the second fundamental form of surface $M_t$ are given by
the equations:
\begin{eqnarray}\label{equation1}
l &=& He^{\lambda}+|z|^n\cos(2t-n~Arg~z),\nonumber\\
m &=& |z|^n\sin(2t-n~Arg~z),\\
n &=& He^{\lambda}-|z|^n\cos(2t-n~Arg~z),\nonumber
\end{eqnarray}
where $t$ is a parameter of the $H$--deformation, cf \cite{Ume}. The following two cases are possible.

\medskip
(i) An $H$--deformation, such that $t$ is constant on $M$.
It can be immediately seen that in new coordinates $\tilde u=\cos t~u+
\sin t~v$, $\tilde v=-\sin t~u+\cos t~v$ the first and the second
forms of surfaces $M_0$ and $M_t$ are the same. By the fundamental
theorem, surfaces $M_0$ and $M_t$ may differ only by a rigid motion
in $E^3$. Thus, the $H$--deformation is trivial.

\medskip
(ii) A non--trivial $H$--deformation. By item (i), $t$ varies
for the points of $M$. Thus far, associated to every $z\in M\backslash 0$, 
there is a chart in which the second fundamental form of surface $M_t(z)$
writes as
\displaymath
l=He^{\lambda}+\cos 2t,\quad
m=\sin 2t,\quad
n=He^{\lambda}-\cos 2t,
\enddisplaymath
where $t$ is the deformation parameter, cf \cite{Wol}.  
A straightforward calculation shows that the principal curvature lines
of the surface $M_t(z)$ coincide with the $\varphi$--geodesic lines of the slope $t$
on $M$. (This fact follows also from the Bonnet Theorem.)
Since every regular point $z\in M$ can be endowed with
such a chart, Lemma \ref{angle} is proved.
$\square$

\section{Proof of Theorem 1}
Take a convex $C^{\infty}$ immersion $f_0:S^2\to E^3$
of the 2--sphere into the Euclidean space $E^3$ which is not totally
umbilic (i.e. there are no $U\subseteq S^2$ such that $f_0(U)$ is
a part of the round sphere). In other words, umbilics are supposed isolated
and their number is finite. Denote by $ds_0$ a Riemann metric on $S^2$
induced by the immersion $f_0$ and by $H:S^2\to\R$ the corresponding mean
curvature function. 
\begin{dfn}
By a Hopf spheroid in $E^3$ we understand a convex $C^{\infty}$ immersion
$f:S^2\to E^3$ such that there exists at least one umbilical point $p$
and a small closed disc $D\ni p$ such that $H(D)=~Const$.
\end{dfn}
\begin{lem}
There exist infinitely many Hopf spheroids in $E^3$. 
\end{lem}
{\it Proof.} By the results of Wente and Kapouleas 
any compact orientable surface $S_g$ of genus $g>0$ admits an immersion into $E^3$
which is a CMC-surface with $H>0$; cf. \cite{Kap}, \cite{Wen}.  
Fix $g\ge 2$ and consider the lines of principal curvature of any such immersion.
By the index argument, there exists an umbilic $p\in S_g$ and a small closed disc $D\ni p$
which is a convex local surface in $E^3$. We separate this local surface from $S_g$.
To obtain a Hopf spheroid, it remains to complete this piece of CMC-surface to a $C^{\infty}$ 
immersion $S^2\to E^3$. By Urysohn's lemma this can be done in an infinite number of ways.
$\square$
\begin{lem}\label{Hopf}
For the Hopf spheroids the Caratheodory conjecture is true.
\end{lem}
{\it Proof.} Without loss of generality we can assume that the umbilic
point $p$ of Hopf spheroid is unique. (For otherwise, if there are
more than one umbilic then we are done.) Since a Hopf spheroid is
locally CMC, we apply  Lemma \ref{angle} to identify the curvature lines in the
disc $D\ni p$ with $\varphi$-geodesic lines in the vicinity of
a singularity $w$.

Let $w\in[w]$ be a word of the minimal length
in the normalized equivalence class $[w]$. According to Lemma
\ref{fol}, there exists a singularity of order $n$ whose
topological type is encoded by the sequence $w$ of symbols $h$ and
$p$. Let $w$ admit $\langle h\rangle$ symbols of type $h$ and
$\langle p\rangle$ symbols of type $p$. By the normalization
axiom, $\langle h\rangle\le n+2$. 

To estimate the Euler--Poincar\'e index of  singularity $w$, note that
the parabolic sectors make no contribution to the index value and
the number $\langle p\rangle$ can be neglected. To the contrary,
if there are no hyperbolic sectors (i.e. $w=p$) we necessarily
have one parabolic sector. The general formula is true:
\displaymath
Ind~w=\cases{1-{\langle h\rangle\over 2}\quad\hbox{if}\quad w\ne p,\cr 
              +1\quad\quad \hbox{if}\quad w=p.}
\enddisplaymath
In either case $Ind~w\le 1$ and by the index argument the conjecture
follows.
$\square$

\bigskip
Now we are ready to finish the proof of Theorem 1. But first we wish to outline the
main idea. To every convex $C^\infty$ immersion $f_0:S^2\to E^3$ one can
relate a Hopf spheroid. This spheroid is uniquely defined by $f_0$
and is  a `modification' of $f_0$ which has an interesting `mechanical' interpretation.

Suppose that $f_0$ is a convex steel ball filled-up with a gas under a pressure. 
Let $p$ be an isolated umbilic of $f_0$. We drill a small hole in $p$ and glue-up
a soap film $D$ into this hole maintaining a pressure
\footnote{The absolute value of the pressure depends on how `flat' is the surface 
at the point $p$. Of course, by  `pressure' we understand difference of pressures
inside and outside the steel ball.}
 inside the ball. We also `deform'
slightly the `edges' of the cut in order to keep the modified surface $f:S^2\to E^3$ in the class
$C^{\infty}$. We claim that $f$ is a Hopf spheroid. 

Indeed, $f(D)$ is a local CMC-surface with an umbilic point $p\in D$. Moreover,
the index of umbilic on the Hopf spheroid is equal to the index of $p$ on $f_0$. 
(This is because the  foliation by  
principal curvature lines at the `steel part' of ball remains intact.) 
In general, if ${\cal F}_0$ and $\cal F$ are foliations by the principal
curvature lines on $f_0$ and $f$, respectively, then ${\cal F}$ is
obtained from ${\cal F}_0$ by a {\it homotopy of openning of a leaf}; cf \cite{N}.

\bigskip
Let $f_0$ be as above. If $p$ is an isolated umbilic of $f_0$ then we take a closed disc
$|D|\le r$ centred at the point $p$. We are going to define a local CMC-surface
$f(D)$. Let $z=u+iv$ be a local parameter which corresponds to a part of CMC-surface
with an umbilic; see the beginning of this section. By the results of Umehara \cite{Ume} 
(see also \cite{Car}, \cite{Che})
 there exists a family of isometric $H$-deformations
depending on a real parameter $t$: 
\begin{equation}\label{equation2}
I=e^{2\lambda}|dz|^2,\qquad II_t=ldu^2+2m~dudv+ndv^2,
\end{equation}
with $l,m$ and $n$ given by equations (\ref{equation1}). 
 The Mainardi-Codazzi and Gauss equations for $I,II_t$:
\begin{equation}
{\partial\varphi\over\overline\partial z} = {\partial H\over\partial z}, \qquad
|\varphi|^2 = e^{4\lambda}(H^2-K),
\end{equation}
where $\varphi =e^{it}z^n$ is a complex quadratic form $\varphi dz^2$,
are satisfied for any real $t$. (Indeed, the first
equation is true since $H=~Const$ and $\varphi$ is holomorphic; the second  
equation follows from $|\varphi e^{it}|=|\varphi|$ and the fact that $H$-deformation
is an isometry.) Therefore, the fundamental forms (\ref{equation2})  are realized
by a concrete local CMC-surface for each real number $t$.

Let $f_t(D)$ be a family of local CMC-surfaces described above. Denote by
$A$ an annular region which surrounds disc $D$:
\begin{equation}\label{equation4}
A=\{z=u+iv|~r\le |z|\le r+\varepsilon\}.
\end{equation}
To glue-up $f_t(D)$ properly, we fix the metric $\lambda$ so that
$\lambda|_{\partial A_{r+\varepsilon}}=\lambda|_{\partial A_r}$,
where the left part denotes a metric on the exterior boundary of $A$ which
is induced by  metric of the surface $f_0$.
The boundary condition $\lambda|_{\partial A_r}$ gives a unique solution
$f_{t=t^*}(D)$ to the Gauss equation, so that a representative in the
family $f_t(D)$ is fixed.

To obtain a $C^{\infty}$ Hopf spheroid it remains to conjugate $f_{t^*}(D)$
with the rest of the sphere:
\begin{equation}\label{equation5}
f(S^2)=\cases{f_{t^*}(D),\quad\hbox{if}
\quad z\in D\subset Int~D_{r+\varepsilon}, \cr  
                         f_0(S^2), \qquad\hbox{if}
\quad z\in S^2\backslash D_{r+\varepsilon}.} 
\end{equation}
By the Urysohn Lemma, function $f$ in formula (\ref{equation5})
can be chosen $C^{\infty}$ for an arbitrary small $\varepsilon$, see
formula (\ref{equation4}). Moreover, taking $r$ sufficiently small
we can fix number $n$ (see (\ref{equation1})) equal to
the order of quadratic form $\varphi$ at point $p$ of the surface $f_0$. 
(Such an order is correctly defined for any $\varphi$, not necessary holomorphic.)

Thus, the surface $f$ given by equation (\ref{equation5}) is a Hopf spheroid.
By the Lawson-Tribuzy theorem $f$ is uniquely defined up to a rigid motion
in $E^3$; see \cite{LaT}. To finish the proof of Caratheodory conjecture,
it remains to notice that  passage from $f_0$ to $f$ gives us a homotopy
$h({\cal F}_0)={\cal F}$ between foliations induced by curvature lines. 
In particular, $Ind~p_0=Ind~p$. By Lemma \ref{Hopf},
the Caratheodory conjecture follows.
$\square$


\enddocument

\begin{figure}\label{Fig}
\begin{center}
\leavevmode            
\epsfxsize=8.cm 
\epsfysize=14.cm
\epsfbox{RISUNKI/jdeq.ps}
\end{center}
\caption{Geodesic foliation near zero of order $1$.}
\end{figure}

\end{document}